\documentclass[12pt,twoside]{amsart}
\usepackage{amsthm,amsfonts,amssymb,amscd,euscript}

\usepackage[matrix,arrow]{xy}

\author{Florin Ambro} 
\address{DPMMS, CMS\\
University of Cambridge,
Wilberforce Road, Cambridge CB3 0WB, UK.}
\email{f.ambro@dpmms.cam.ac.uk}

\newcommand{\Q}{{\mathbb Q}}
\newcommand{\Z}{{\mathbb Z}}
\newcommand{\N}{{\mathbb N}}
\newcommand{\R}{{\mathbb R}}

\newcommand{\calE}{{\mathcal E}}

\newcommand{\calO}{{\mathcal O}}

\newcommand{\bA}{{\mathbf A}}

\newcommand{\bD}{{\mathbf D}}
\newcommand{\bE}{{\mathbf E}}

\newcommand{\bH}{{\mathbf H}}
\newcommand{\bK}{{\mathbf K}}

\newcommand{\bM}{{\mathbf M}}

\newcommand{\bP}{{\mathbf P}}

\newcommand{\Exc}{\operatorname{Exc}}

\newcommand{\mult}{\operatorname{mult}}

\newcommand{\Pic}{\operatorname{Pic}}

\newcommand{\Proj}{\operatorname{Proj}}

\newcommand{\Spec}{\operatorname{Spec}}
\newcommand{\Supp}{\operatorname{Supp}}
\newcommand{\Var}{\operatorname{Var}}

\theoremstyle{plain}
\newtheorem{thm}{Theorem}[section]

\newtheorem{lem}[thm]{Lemma}

\newtheorem{prop}[thm]{Proposition}

\newtheorem{conj}[thm]{Conjecture}

\theoremstyle{definition}
\newtheorem{defn}[thm]{Definition}

\newtheorem{example}{Example}

\newtheorem{rem}[thm]{Remark}

\newtheorem{ack}{Acknowledgments}  
  
\newtheorem{ABC}{Abundance Conjecture}

\theoremstyle{remark}

\begin{document}

\bibliographystyle{amsalpha+}
\title[Nef dimension of minimal models]
{Nef dimension of minimal models}

\begin{abstract} We reduce the Abundance
Conjecture in dimension $4$ to the following
numerical statement: if the canonical divisor
$K$ is nef and has maximal nef dimension, 
then $K$ is big. From this point of view, we 
``classify'' in dimension $2$ nef divisors 
which have maximal nef dimension, but which 
are not big.
\end{abstract}

\maketitle


\setcounter{section}{-1}


\section{Introduction}


\footnotetext[1]{1991 Mathematics Subject 
Classification. Primary: 14E30. Secondary: 14J10, 
14N30.}

A {\em minimal model} is a complex projective 
variety $X$ with at most terminal singularities, 
whose canonical divisor $K$ is numerically effective 
(nef): $K\cdot C\ge 0$ for every curve $C\subset 
X$. Up to dimension three, minimal models have a 
geometrical characterization (Kawamata~\cite{ab,Ka92}, 
Miyaoka~\cite{Miy1, Miy2, Miy3}) :

\begin{ABC}\cite{KMM} 
Let $X$ be a minimal model. Then the 
linear system $|kK|$ is base point free, for some 
positive integer $k$. 
\end{ABC}

In dimension four, it is enough to show that $X$ 
has positive Kodaira dimension if $K$ is not 
numerically trivial (Kawamata~\cite{ab}, 
Mori~\cite{Mo}). 

A direct approach is to first construct the 
morphism associated to the expected base point 
free pluricanonical linear systems:
$$
f\colon X \to \Proj(\oplus_{k\ge 0} H^0(X,kK)).
$$
Since $K$ is nef, $f$ is the unique morphism
with connected fibers which contracts exactly
the curves $C \subset X$ with $K\cdot C=0$.
Tsuji~\cite{T00} and Bauer et al~\cite{nr} have 
recently solved this existence problem 
{\em birationally}: 
for any nef divisor $D$ on $X$, there exists a 
rational dominant map $f\colon X-\to Y$ such that 
$f$ is regular over the generic point of $Y$ and 
a very general curve $C$ is contracted by $f$ if 
and only if $D\cdot C=0$. This rational map is 
called the {\em nef reduction} of $D$, and 
$n(X,D):=\dim(Y)$ is called the 
{\em nef dimension} of $D$. 
The nef reduction map is non-trivial, except for
the two extremal cases:
\begin{enumerate}
\item[(i)] $n(X,K)=0$: $K$ 
is numerically trivial in this case~\cite{nr}, 
and Abundance is known (Kawamata~\cite{minmod}).
\item[(ii)] $n(X,K)=\dim(X)$:
$K\cdot C>0$ for very general curves $C \subset
X$.
The nef reduction rational map is the identity.
\end{enumerate}

Our main result is that Abundance holds for a 
minimal model $X$ if the nef reduction map is 
non-trivial and the Log Minimal Model Program 
and Log Abundance hold in dimension $n(X,K)$. The
latter two conjectures are known to hold up to 
dimension three (Shokurov~\cite{Sh}, Keel, Matsuki, 
McKernan~\cite{LAB}), hence we obtain

\begin{thm}\label{main} 
Let $X$ be a minimal model with $n(X,K)\le 3$.
Then the linear system $|kK|$ is base point free 
for some positive integer $k$.
\end{thm}

The Base Point Free Theorem 
(Kawamata, Shokurov~\cite{KMM}) states that
Abundance holds if the canonical class $K$ is big.
Combined with Theorem~\ref{main}, the $4$-dimensional
case of Abundance is equivalent to the following 

\begin{conj}\label{q} Let $X$ be a minimal
$4$-fold. If $K$ has maximal nef dimension, 
then $K$ is big.
\end{conj}

We stress that this statement is {\em numerical}:
since $K$ is nef, $K$ is big if and only if
$K^{\dim(X)}\ne 0$. For this reason, it is 
important to investigate how far are 
(adjoint) divisors of maximal nef dimension 
from being big. 
Questions of similar type have appeared in the
literature: a divisor $D$ is {\em strictly 
nef} (Serrano~\cite{Sr1}) if $D\cdot C>0$ for every 
curve $C\subset X$. Up to dimension $3$, it is known 
that $\pm K$ is strictly nef if and only if $\pm K$ 
is ample (see ~\cite{Sr1,UH} and the references 
there). We point out that Conjecture~\ref{q} is 
false for the anti-canonical divisor $-K$ (which,
at least in dimension two, is the only exception
below):

\begin{thm}\label{sf} 
Let $X$ be a smooth projective surface.
Assume that $D$ is a nef Cartier divisor of maximal 
nef dimension, which is not big. 
Then exactly one of the following cases occurs:
\begin{itemize}
\item[(1)] The divisor $K+tD$ is big for $t>2$.
\item[(2)] There exists a birational contraction 
$f\colon X\to Y$ and there exists $t \in (0,2]$ such 
that $D=f^*(D_Y)$, and $K_Y+tD_Y\equiv 0$. Moreover,
$D$ is effective up to algebraic equivalence.
In Sakai's classification table~\cite{Sakai83}, $Y$ is 
either a degenerate Del Pezzo, or an elliptic rulled 
surface of type $II_c,II_c^*$.
\end{itemize}
\end{thm}

Theorem~\ref{main} is proved in several steps. The 
properties of the nef reduction map $f$ and the 
numerically trivial case of Abundance~\cite{minmod} 
imply that $f$ is birational to a parabolic fiber 
space $f'\colon X'\to Y'$, and the canonical class 
$K$ descends to a divisor $P$ on $Y'$. 
After an idea of Fujita~\cite{Fuj86}, it is enough 
to show that $P$ is the semi-positive part in the 
Fujita decomposition associated to a 
log variety $(Y',\Delta)$: the semi-ampleness of
$P$ follows then from the Log Minimal Model Program
and Log Abundance applied to $(Y',\Delta)$.
The key ingredient in this argument is an adjunction 
formula for the parabolic fiber space $f'$
(Kawamata~\cite{Ka83, minmod}, Fujino, Mori~\cite{fm,
osamu}), similar to Kodaira's formula for elliptic 
surfaces. We expect that the logarithmic version of 
Theorem~\ref{main} follows from the same argument,
provided that Kawamata's adjunction 
formula~\cite{Ka83} is extended to the logarithmic 
case (see also Fukuda~\cite{Fuk}).

Finally, Theorem~\ref{sf} follows from the 
classification of surfaces and generalizes 
a result of Serrano~\cite{Sr1}.

\begin{ack} This work was supported through a 
European Community Marie Curie Fellowship.
\end{ack}


\section{Preliminary}


A {\em variety} is a reduced and irreducible 
separable scheme of finite type, defined over an 
algebraically closed field of characteristic zero. 
A {\em contraction} is a proper morphism $f\colon X
\to Y$ such that $\calO_Y=f_*\calO_X$.

Let $X$ a normal variety, and let $K\in \{\Z,\Q,\R\}$. 
A $K$-Weil divisor is an element of $Z^1(X)\otimes_\Z K$.
Two $\R$-Weil divisors $D_1, D_2$ are $K$-{\em linearly 
equivalent}, denoted $D_1\sim_K D_2$, if there exist 
$q_i\in K$ and rational functions $\varphi_i\in k(X)^\times$ 
such that $D_1-D_2=\sum_i q_i(\varphi_i)$. An $\R$-Weil 
divisor $D$ is called
\begin{itemize}
\item[(i)] {\em $K$-Cartier} if $D\sim_K 0$ in 
a neighborhood of each point of $X$.
\item[(ii)] {\em nef} if $D$ is $\R$-Cartier and 
$D\cdot C\ge 0$ for every curve $C\subset X$.
\item[(iii)] {\em ample} if $X$ is projective and 
the numerical class of $D$ belongs to the real cone 
generated by the numerical classes of ample Cartier 
divisors.
\item[(iv)] {\em semi-ample} if there exists a 
contraction $\Phi\colon X\to Y$ and an ample 
$\R$-divisor $H$ on $Y$ such that $D\sim_\R \Phi^*H$. 
If $D$ is rational, this is equivalent to the
linear system $|kD|$ being base point free for some $k$.
\item[(v)] {\em big} if there exists $C>0$ such that
$\dim H^0(X,kD)\ge Ck^{\dim(X)}$ for $k$ sufficiently 
large and divisible. By definition,
$$
H^0(X,kD)=\{a\in k(X)^\times; (a)+kD\ge 0\}\cup \{0\}. 
$$
\end{itemize}

The {\em Iitaka dimension} of $D$ is 
$
\kappa(X,D)=\max_{k\ge 1} \dim \Phi_{|kD|}(X),
$
where $\Phi_{|kD|}\colon X-\to {\mathbb P}(|kD|)$ 
is the rational map associated to the linear system 
$|kD|$. If all the linear systems $|kD|$ are empty, 
$\kappa(X,D)=-\infty$. If $D$ is nef, the 
{\em numerical dimension} $\nu(X,D)$ is the largest 
non-negative integer $k$ such that there exists a 
codimension $k$ cycle $C\subset X$ with 
$D^k\cdot C=0$. 

\begin{defn} (V.V. Shokurov) A {\em $K$-b-divisor} 
$\bD$ of $X$ is a family $\{\bD_{X'}\}_{X'}$ of 
$K$-Weil divisors indexed by all birational models  
of $X$, such that $\mu_*(\bD_{X''})=\bD_{X'}$ if
$\mu\colon X''\to X'$ is a birational contraction. 

Equivalently, $\bD=\sum_E \mult_E(\bD) E$ is a $K$-valued 
function on the set of all (geometric) valuations of the 
field of rational functions $k(X)$, having finite support 
on some (hence any) birational model of $X$.
\end{defn}

\begin{example} (1) Let $\omega$ be a top rational 
differential form of $X$. The associated family of divisors 
$\bK=\{(\omega)_{X'}\}_{X'}$ is called the {\em canonical 
b-divisor} of $X$. 

(2) A rational function $\varphi \in k(X)^\times$ defines a 
b-divisor $\overline{(\varphi)}=\{(\varphi)_{X'}\}_{X'}$.

(3) An $\R$-Cartier divisor $D$ on a birational model $X'$ 
of $X$ defines an $\R$-b-divisor $\overline{D}$ such that
$(\overline{D})_{X''}=\mu^*D$ for every birational contraction 
$\mu\colon X''\to X'$. 
\end{example}

An $\R$-b-divisor $\bD$ is called {\em $K$-b-Cartier} 
if there exists a birational model $X'$ of $X$ such that
$\bD_{X'}$ is $K$-Cartier and $\bD=\overline{\bD_{X'}}$.
In this case, we say that $\bD$ {\em descends to $X'$}.
An $\R$-b-divisor $\bD$ is {\em b-nef}
({\em b-semi-ample}, {\em b-big}, {\em b-nef and good}) 
if there exists a birational contraction $X'\to X$ such 
that $\bD=\overline{\bD_{X'}}$, and $\bD_{X'}$ is 
nef (semi-ample, big, nef and good).

A {\em log pair} $(X,B)$ is a normal variety $X$ 
endowed with a $\Q$-Weil divisor $B$ such that $K+B$ 
is $\Q$-Cartier. A {\em log variety} is a log pair 
$(X,B)$ such that $B$ is effective. The 
{\em discrepancy $\Q$-b-divisor} of a log pair $(X,B)$ 
is
$$
\bA(X,B)=\bK-\overline{K+B}.
$$ 
A log pair $(X,B)$ is said to have at most 
{\em Kawamata log terminal singularities} if 
$\mult_E(\bA(X,B))> -1$ for every geometric 
valuation $E$.


\section{Nef reduction}


The existence of the nef reduction map is 
originally due to Tsuji~\cite{T00}. An algebraic 
proof of the sharper statement below is due to Bauer, 
Campana, Eckl, Kebekus, Peternell, Rams, Szemberg, 
and Wotzlaw~\cite{nr}.

\begin{thm}\cite{T00,nr} 
Let $D$ be a nef $\R$-Cartier divisor on a normal
projective variety $X$. Then there exists a rational 
map $f \colon X -\to Y$ to a normal projective 
variety $Y$, satisfying the following properties:
\begin{itemize}
\item[(i)]
$f$ is a dominant rational map with connected fibers,
which is a morphism over the general point of $Y$.
\item[(ii)] There exists a countable intersection $U$ of 
Zariski open dense subsets of $X$ such that for every 
curve $C$ with $C\cap U\ne \emptyset$, $f(C)$ is a point
if and only if $D\cdot C=0$.
\end{itemize}
In particular, $D|_W\equiv 0$ for general fibers $W$
of $f$.
\end{thm}

 The rational map $f$ is unique, and is called the 
{\em nef reduction of } $D$. The dimension of $Y$ is 
called the {\em nef dimension} of $D$, denoted by $n(X,D)$. 
In general, the following inequalities hold~\cite{ab,nr}:
$$
\kappa(X,D)\le \nu(X,D)\le n(X,D)\le \dim(X).
$$

\begin{defn} A nef $\Q$-Cartier divisor $D$ is 
called good if 
$$
\kappa(X,D)=\nu(X,D)=n(X,D).
$$
\end{defn}

\begin{rem} This is equivalent to Kawamata's 
definition~\cite{ab}. If $\kappa(X,D)=\nu(X,D)$,
there exists a dominant rational map $f\colon 
X-\to Y$ and a nef and big $\Q$-divisor $H$ on $Y$ 
such that $\overline{D}\sim_\Q f^*(\overline{H})$,
by ~\cite{ab}. Thus $n(X,D)$ coincides with the
Iitaka and numerical dimension. 
\end{rem}

\begin{rem}\cite{nr} The extremal values 
of the nef dimension are:
\begin{itemize}
\item[(i)] $n(X,D)=0$ if and only if $D$ is 
numerically trivial ($\nu(X,D)=0$).
\item[(ii)] $n(X,D)=\dim(X)$ if and only if 
there exists a countable intersection $U$ of Zariski
open dense subsets of $X$ such that $D\cdot C>0$ for
every curve $C$ with $C\cap U\ne \emptyset$.
\end{itemize}
\end{rem}


\section{Fujita decomposition}


\begin{defn}\cite{Fuj86} 
An $\R$-Cartier divisor $D$ on a normal proper 
variety $X$ has a {\em Fujita decomposition} if 
there exists a b-nef $\R$-b-divisor $\bP$ of $X$
with the following properties:
\begin{itemize}
\item[(i)] $\bP\le \overline{D}$.
\item[(ii)] $\bP=\sup\{\bH; \bH \mbox{ b-nef 
$\R$-b-divisor}, \bH \le \overline{D}\}$.
\end{itemize}
The $\R$-b-divisor $\bP=\bP(D)$ is unique if 
it exists, and is called the {\em semi-positive 
part} of $D$. The $\R$-b-divisor 
$\bE=\overline{D}-\bP$ is called the 
{\em negative part} of $D$, and
$
\overline{D}=\bP+\bE
$
is called the {\em Fujita decomposition} of
$D$.
\end{defn}

\begin{rem}
Allowing divisors with real coefficients is
necessary: there exist Cartier divisors
(in dimension at least $3$) which have a Fujita
decomposition with irrational semi-positive 
part~\cite{Cut}.
\end{rem}

Clearly, a nef $\R$-Cartier divisor $D$ has a 
Fujita decomposition, with semi-positive 
part $\overline{D}$. More examples can be 
constructed using the following property:

\begin{prop}\cite{Fuj86}\label{in} 
Let $f\colon X\to Y$ be a proper contraction, let
$D$ be an $\R$-Cartier divisor on $Y$ and let $E$ 
be an effective $\R$-Cartier divisor on $X$ such that $E$ 
is vertical and supports no fibers over codimension one 
points of $Y$. 

Then $D$ has a Fujita decomposition if and only if 
$f^*D+E$ has a Fujita decomposition, and moreover, 
$\bP(f^*D+E)=f^*(\bP(D))$.
\end{prop}

\begin{lem} Assume LMMP and Log Abundance. Let 
$(X,B)$ be a log variety with log canonical 
singularities. Then $K+B$ has a Fujita decomposition 
if and only if $\kappa(X,K+B)\ge 0$, and the 
semi-positive part is semi-ample. Moreover, 
$$
\bP(K+B)=\overline{K_Y+B_Y},
$$
for a log minimal model $(Y,B_Y)$.
\end{lem}

\begin{proof} If $K+B$ is nef, it has a Fujita 
decomposition with semi-positive part 
$\overline{K+B}$. By Abundance, it is semi-ample.
If $K+B$ is not nef, we run the LMMP for $(X,B)$. 
We may assume that $X$ is $\Q$-factorial by 
Proposition~\ref{in}. If $f\colon (X,B)\to Y$ is a 
divisorial contraction, then 
$$
K+B=f^*(K_Y+B_Y)+\alpha E,
$$
where $E$ is exceptional on $Y$ and $\alpha>0$. 
Thus $K+B$ has a Fujita decomposition if and only 
if $K_Y+B_Y$ has, and the semi-positive parts 
coincide. If $t\colon (X,B)-\to (X^+,B_{X^+})$ 
is a log-flip,
$$
\overline{K+B}=\overline{K_{X^+}+B_{X^+}}+\bE,
$$
where $\bE$ is an effective $\Q$-b-divisor which
is exceptional on both $X$ and $X'$. Therefore $K+B$ 
has a Fujita decomposition if and only if 
$K_{X^+}+B_{X^+}$ has, and the semi-positive parts 
coincide.

If $f\colon (X,B)\to Y$ is a log Fano fiber space,
$K+B$ admits no Fujita decomposition. 
\end{proof}

\begin{lem}\label{vn} \cite{ab,Fuj86}  
Let $f\colon X\to Y$ be a contraction of normal 
proper varieties, and let $D$ be a nef $\R$-divisor
on $X$ which is vertical on $Y$. Then there exists a
b-nef $\R$-b-divisor $\bD$ of $Y$ such that
$\overline{D}=f^*\bD$.
\end{lem}

\begin{proof} After a resolution of singularities, 
Hironaka's flattenining and the normalization of the 
total space of the induce fibration, we have a fiber 
space induced by birational base change
\[ 
\xymatrix{
X\ar[d]_f    & X'\ar[d]_{f'} \ar[l]_\mu     \\
Y     &         Y'\ar[l]
}\]
such that $f'$ is equi-dimensional, $X'$ is 
normal and $Y'$ is non-singular, and $\mu^*D$ is 
vertical on $Y'$. Let $D'$ be the largest $\R$-divisor
on $Y'$ such that ${f'}^*D'\le \mu^*D$. Since
$f'$ is equi-dimensional, $E=\mu^*D-{f'}^*D'$ is 
effective and supports no fibers over codimension 
one points of $Y'$.
Furthermore, $E$ is $f'$-nef since $D$ is nef.
By \cite[Lemma 1.5]{Fuj86}, $E=0$. 
Therefore $\mu^*D={f'}^*(D')$. In particular, 
$D'$ is nef and $\bD=\overline{D'}$ satisfies
the required properties.
\end{proof}


\section{Parabolic fiber spaces}


 We recall results of Kawamata~\cite{Ka83,minmod} 
and Fujino, Mori~\cite{fm,osamu} on adjunction 
formulas of Kodaira type for parabolic fiber spaces. 
Their results are best expressed through Shokurov's 
terminology of b-divisors. With a view towards the
logarithmic case, we introduce them via lc-trivial
fibrations (see~\cite{bp}).

A {\em parabolic fiber space} is a contraction
of non-singular proper varieties $f\colon X\to Y$ 
such that the generic fiber $F$ has Kodaira 
dimension zero. Let $b$ be the smallest positive 
integer with $|bK_F|\ne \emptyset$. We fix a rational 
function $\varphi \in k(X)^\times$ such that 
$K+\frac{1}{b}(\varphi)$ is effective over the 
generic point of $Y$.

\begin{lem} 
There exists a unique $\Q$-divisor $B_X$ on $X$ 
satisfying the following properties:
\begin{itemize}
\item[(i)] $K_X+B_X+\frac{1}{b}(\varphi)=f^*D$ 
for some $\Q$-divisor $D$ on $Y$.
\item[(ii)] There exists a big open subset 
$Y^\dagger \subseteq Y$ such that 
$-B_X|_{f^{-1}(Y^\dagger)}$ is effective and contains 
no fibers of $f$ in its support.
\end{itemize}
In particular, $f\colon (X,B_X)\to Y$ is an lc-trivial
fibration.
\end{lem}

\begin{defn} Let $f\colon X\to Y$ be a parabolic
fiber space with a choice of a rational function 
$\varphi$, as above. The {\em moduli $\Q$-b-divisor} 
of $f$, denoted $\bM=\bM(f,\varphi)$, is the moduli 
$\Q$-b-divisor of the lc-trivial fibration $f\colon 
(X,B_X)\to Y$.
\end{defn}

If $\varphi'$ is another choice of the rational 
function, then $b\bM(f,\varphi)\sim b\bM(f,\varphi')$. 
Therefore $b\bM$ is uniquely defined up to linear 
equivalence. According to the following Lemma, $\bM$ 
is independent of birational changes of $f$:

\begin{lem} Consider a commutative diagram
\[ 
\xymatrix{
 X \ar[d]_f  & X' \ar[l]_\nu \ar[d]_{f'}   \\
 Y  & Y' \ar[l]_\mu
} \]
where $f,f'$ are parabolic fiber spaces and $\mu,\nu$
are birational contractions. Then $\bM(f)=\bM(f')$.
\end{lem}

\begin{proof} Assume first that $\mu$ is the identity
morphism. Since $X,X'$ are nonsingular, it is easy to 
see that $\bA(X,B_X)=\bA(X',B_{X'})$. 
Therefore $\bM(f)=\bM(f')$.

We are left with the case when $\nu$ is the 
identity morphism. Let $B_X^{(Y)}$ and $B_X^{(Y')}$ 
be the $\Q$-divisors induced by $f$ and $f'$, 
respectively. 
Since the general fibre is non-singular of zero Kodaira 
dimension, there exists a $\Q$-divisor $C$ on $Y'$ such 
that $B_X^{(Y')}=B_X^{(Y)}+{f'}^*C$. 
Therefore $\bM(f)=\bM(f')$, by ~\cite[Remark 3.3]{bp}.
\end{proof}

\begin{prop}\label{pb}
Let $f\colon X\to Y$ be a parabolic fiber space. 
\begin{itemize}
\item[(1)] Consider a commutative diagram
\[ 
\xymatrix{
 X \ar[d]_f  & X' \ar[l]_\nu \ar[d]_{f'}   \\
 Y  & Y' \ar[l]_\varrho
} \]
where  $\varrho$ is a surjective proper morphism, and 
$f'$ is an induced parabolic fiber space. Then 
$
\varrho^*\bM(f)\sim_\Q \bM(f').
$

\item[(2)] If $f$ is semi-stable in codimension 
one, then
$$
f_*\calO_X(iK_{X/Y})^{**}=\calO_Y(i\bM_Y) \cdot 
\varphi^i, \mbox{ for } b|i 
$$
\item[(3)] The moduli $\Q$-b-divisor $\bM(f)$ is
b-nef.
\end{itemize}
\end{prop}

The key result of this section is the following
corollary of~\cite[Theorem 3.6]{minmod}:

\begin{thm}\label{k0}
Let $f\colon X\to Y$ be a parabolic fiber space. 
Assume that its geometric generic fibre 
$X\times_Y \Spec(\overline{k(Y)})$ is birational to a 
normal variety $\bar{F}$ with canonical singularities, 
defined over $\overline{k(Y)}$, such that $K_{\bar{F}}$ 
is semi-ample.
Then the moduli $\Q$-b-divisor $\bM(f)$ is b-nef and
good.
\end{thm}

\begin{proof} From the definiton of the variation of a 
fibre space, there exists a commutative diagram

\[ \xymatrix{
X \ar[d]_f & \bar{X}\ar[l]\ar[d]_{\bar{f}}
\ar[r] & X^!\ar[d]_{f^!}\\
Y        & \bar{Y} \ar[l]_\tau \ar[r]^\varrho & Y^!
} \]

such that the folowing hold:
\begin{itemize}
\item[(1)] $\bar{f}$ and $f^!$ are parabolic fiber spaces.
\item[(2)] $\tau$ is generically finite, and $\varrho$
is a proper dominant morphism.
\item[(3)] $\bar{f}$ is birationally induced via base 
change by both $f$ and $f^!$.
\item[(4)] $\Var(f)=\Var(f^!)=\dim(Y^!)$.
\end{itemize}

Let $\bM,\bar{\bM},\bM^!$ be the corresponding moduli 
$\Q$-b-divisors. After a generically finite base change, 
we may also assume that $\bM^!$ descends to $Y^!$, and 
$f^!$ is semi-stable in codimension one.
By (3) and Proposition~\ref{pb}, we have
$$
\tau^*\bM=\bar{\bM}\sim_\Q \varrho^*(\bM^!).
$$
In particular, $\kappa(\bM)=\kappa(\bM^!)$. Since 
$\bar{F}$ is a good minimal model, Viehweg's $Q(f^!)$ 
Conjecture holds~\cite[Theorem 1.1.(i)]{minmod}, that is 
the sheaf
$
(f^!_*\omega_{X^!/Y^!}^i)^{**}
$
is big for $i$ large and divisible. But 
$
(f^!_*\omega_{X^!/Y^!}^i)^{**} \simeq \calO_{Y^!}
(i\bM^!_{Y^!})
$
for $b|i$, since $f^!$ is semi-stable in codimension 
one. Equivalently, $\kappa(Y^!,\bM^!_{Y^!})=
\dim(Y^!)$, or $\bM^!$ is b-nef and big.
Therefore $\tau^*\bM$ is b-nef and good, hence $\bM$ 
is b-nef and good.
\end{proof}


\section{Reduction argument}


\begin{thm} Let $X$ be a projective variety with 
canonical singularities such that the canonical
divisor $K$ is nef. If $n(X,K)\le 3$, then the 
canonical divisor $K$ is semi-ample.
\end{thm}

\begin{proof} Let $\Phi\colon X -\to Y$ be the
quasi-fibration associated to the nef canonical 
divisor $K$ of $X$, and let $\Gamma$ be the 
normalization of the graph of $\Phi$:
\[ \xymatrix{
 & \Gamma \ar[dl]_\mu \ar[dr]^f & \\
  X      & & Y
} \] 
Since $\Phi$ is a quasi-fibration, $\mu$ is 
birational, $f$ is a contraction and $\Exc(\mu)
\subset \Gamma$ is vertical over $Y$. Let
$W$ be a general fibre of $f$.

{\em Step 1}: $W$ is a normal variety with canonical 
singularities, and $K_W\sim_\Q 0$.
Indeed, $W$ has canonical singularities and 
$K_W= \mu^*K|_W$. 
The definition of $\Phi$ implies that $K_W$ is 
numerically trivial. From ~\cite[Theorem 8.2]{minmod}, 
we conclude that $K_W\sim_\Q 0$. 

{\em Step 2}: There exist a diagram
\[ \xymatrix{
 X & X' \ar[l]_\mu \ar[d]^{f'}  \\
   &  Y'
} \]
satisfying the following properties:
\begin{enumerate}
\item[(a)] $\mu$ is a birational contraction.
\item[(b)] $f'\colon X'\to Y'$ is a parabolic
fiber space.
\item[(c)] There exists a simple normal crossings
divisor $\Sigma$ on $Y'$ such that $f'$ is smooth
over $Y'\setminus \Sigma$.
\item[(d)] The moduli $\Q$-b-divisor $\bM=\bM(f')$ 
descends to $Y'$ and there exists a contraction 
$h\colon Y'\to Z$ and a nef and big $\Q$-divisor 
$N$ on $Z$ such that $\bM_{Y'}\sim_\Q h^*N$.
\item[(e)] Let $E$ be any prime divisor on $X'$.
If $E$ is exceptional over $Y'$, then $E$ is 
exceptional over $X$.
\end{enumerate}

Indeed, we may assume that $Y$ is non-singular. Let 
$\Gamma'\to \Gamma$ be a resolution of singularities, 
and let $f_0\colon \Gamma' \to Y$ be the induced 
contraction. The general fiber of $f_0$ is birational
to the general fiber of $f$. The latter is a normal
variety $W$ with canonical singularities, and 
$K_W\sim_\Q 0$. Therefore $f_0$ is a parabolic 
fiber space. We define $f'\colon X'\to Y'$ to be
a parabolic fiber space induced after a sufficiently
large birational base change $Y'\to Y$. 
By Theorem~\ref{k0}, the moduli $\Q$-b-divisor 
$\bM(f)=\bM(f_0)$ satisfies (d) once $Y'$ dominates
a certain resolution of $Y$. Also, (e) holds once 
$f'$ dominates a flattening of $f$, and (b) follows
from Hironaka's embedded resolution of singularities. 

{\em Step 3}: 
There exists a effective $\Q$-divisor $\Delta$ on 
$Y'$ such that $(Y',\Delta)$ is a log variety with
Kawamata log terminal singularities, $K_{Y'}+\Delta$
has a Fujita decomposition and $\overline{K}\sim_\Q 
{f'}^*(\bP(K_{Y'}+\Delta))$.

Indeed, the parabolic fiber space $f'$ induces an
lc-trivial fibration $(X',B_{X'})\to Y'$, with 
associated discriminant divisor $B_{Y'}$. We have
$$
K_{X'}+B_{X'}+\frac{1}{b}(\varphi)={f'}^*
(K_{Y'}+B_{Y'}+\bM_{Y'}).
$$
It is clear that $B_{Y'}$ is effective,
$\lfloor B_{Y'} \rfloor =0$ and $\Supp(B_{Y'})
\subseteq \Sigma$.
Therefore $(Y',B_{Y'})$ is a log variety 
with Kawamata log terminal singularities.
By (d), there exists an effective $\Q$-divisor 
$\Delta$ on $Y'$ such that $(Y',\Delta)$ is a log 
variety with Kawamata log terminal singularities, 
and $\Delta\sim_\Q B_{Y'}+\bM_{Y'}$. In particular,
$$
K_{X'}+B_{X'} \sim_\Q {f'}^*(K_{Y'}+\Delta).
$$
Let $\mu^*K=K_{X'}-A$ and let $B_{X'}=E^+-E^{-}$ be
the decomposition into positive and negative parts.
It is clear that $A$ is effective and exceptional 
over $X$, and $A-E^-$ is vertical on $Y$. Thus
there exist effective $\Q$-divisors $A'\le A$ and
$E'\le E^-$ such that $A-E^-=A'-E'$ and $E'$ is
vertical and supports no fibers over codimension
one points of $Y'$. In particular,
$$
\mu^*K+A'+E^+\sim_\Q {f'}^*(K_{Y'}+\Delta)+E'.
$$
By (e), the left hand side has a Fujita 
decomposition, with semi-positive part $\overline{K}$. 
Proposition~\ref{in} applies, hence $K_{Y'}+\Delta$ has
a Fujita decomposition and  
and $\overline{K}\sim_\Q  {f'}^*(\bP(K_{Y'}+\Delta))$.

{\em Step 4}: 
From the LMMP and Abundance applied to the log variety
$(Y',\Delta)$, the semi-positive part of $K_{Y'}+\Delta$ 
is b-semi-ample. Therefore $\overline{K}$ is b-semi-ample, 
that is $K$ is a semi-ample $\Q$-divisor.
\end{proof}


\section{Maximal nef dimension which are not big}


We prove Theorem~\ref{sf} in this section. We fix
the notation: $X$ is a smooth projective surface 
and $D$ is a nef Cartier divisor which has
maximal nef dimension, but it is not big. We denote
by $K$ the canonical divisor of $X$.

\begin{prop} \label{srk}
The following hold:
\begin{itemize}
\item[(1)] $\kappa(X,D)\le 0$, $\nu(X,D)=1$.
\item[(2)] $D\cdot K\ge 0$.
\item[(3)] If $D\cdot K=0$, one of the following 
holds:
 \begin{itemize}
\item[a)] $\kappa(X,D)=-\infty$ and $X$ is birational 
to ${\mathbb P}_C(\calE)$, where $C$ is a non-rational
curve.
\item[b)] $\kappa(X,D)=0$ and $X$ is either a rational
surface, or an elliptic ruled surface.
 \end{itemize}
\item[(4)] Assume $D\cdot K=0$ and $K^2\ge 0$. Then $D$ 
is algebraically equivalent to an effective divisor.
\end{itemize}
\end{prop}

\begin{proof} Since $D$ cannot be good, (1)
holds. We have 
$$
\chi(X,mD)=\frac{-D\cdot K}{2}m +\chi(\calO_X).
$$ 
Since $\nu(X,D)>0$, $h^2(mD)=h^0(K-mD)=0$ for $m\gg 0$. 

(2) If $D\cdot K<0$, then $\kappa(X,D)\ge 1$. This
contradicts (1). 

(3) Assume $D\cdot K=0$. In particular, $\kappa(X)\le 0$.
Indeed, let $L$ be a divisor such that $DL=0$. Since $D$ 
is nef, $D$ is orthogonal on the irreducible components of 
all divisors in $|mL|, \ m\ge 0$. Since $D$ is orthogonal 
on at most a countable number of curves, $\kappa(X,L)\le 0$.

Assume $\kappa(X)=0$. Let 
$\sigma\colon X\to X'$ be the birational contraction to 
a minimal model. Since $K_{X'}\sim_\Q 0$, $K\sim_\Q E$
where $E$ is effective and $\Supp(E)=\Exc(\sigma)$. 
Since $D\cdot K=0$, $D$ is orthogonal on each exceptional
divisor, hence $D=\sigma^*(D_{X'})$.
Thus we may assume $X$ is a minimal model. After an 
\'etale cover, $X$ is an Abelian surface or a $K3$ surface.
If $X$ is an Abelian surface, $D$ is big by the same 
argument as in ~\cite[Proposition 1.4]{Sr1}. Contradiction.
If $X$ is a $K3$ surface, $h^0(X,mD)=h^1(X,mD)+2$ by 
Riemann-Roch, hence $\kappa(X,D)\ge 1$. Contradiction.

Therefore $\kappa(X)=-\infty$. Riemann-Roch gives
$$
h^0(X,mD)=h^1(X,mD)+1-q(X), \ \ m \ge 1
$$
If $q(X)=0$, then $h^0(D)>0$. We are in case (b), 
and the rest of the claim is well known (see~\cite{Sakai85}).
Assume $q(X)>0$. Then there exists a birational 
contraction $X\to X'={\mathbb P}_C(E)$, with $q(X)=g(C)
\ge 1$. We are in case (a).

(4) If $q(X)=0$, $|D|\ne \emptyset$ by Riemann-Roch. 
Assume $q(X)>0$. There exists a birational contraction
$X\to X'={\mathbb P}_C(E)$, with $q(X)=g(C)$. Since
$0\le K_X^2\le K_{X'}^2=8(1-q(X))\le 0$, we infer that
$X={\mathbb P}_C(\calE)$ and $C$ is an elliptic curve, 
i.e. $q(X)=1$.

If $h^1(D+F_t-F_0)>0$ for some $t\in C$, then
$h^0(D+F_t-F_0)>0$ by Riemann-Roch. Assume
$h^1(D+F_t-F_0)=0$ for every $t\in C$. Since
$D$ is of maximal nef dimension, $D\cdot F_0>0$.
Therefore $h^0(F_0,D|_{F_0})>0$. By~\cite{Sr1}, 
Proposition 1.5, $D$ is algebraically equivalent to
an effective divisor.
\end{proof}

\begin{thm}\label{sk}
\cite{Sakai85} In the case (3b) above, assume moreover
that $D$ is effective and $DC>0$ for every $(-1)$-curve 
$C$ of $X$. Then the pair $(X,D)$ is classified as 
follows:
\begin{itemize}
\item[(i)] $X$ is a rational surface
such that $-K$ is nef and $K^2=0$:
$$
\kappa(X,-K)=0, \nu(X,-K)=1, n(X,-K)=2. 
$$
There exists a connected effective cycle 
$\sum n_i C_i \in |-K|$ such that the greatest common
divisor of the $n_i$'s is $1$. 
Also, $D=m\sum n_i C_i$ for some positive integer $m$.

\item[(ii)] $X={\mathbb P}_C(\calE)$ is a geometrically
ruled surface over an elliptic curve $C$, of type
$II_c$ or $II_c^*$ in Sakai's classification table:
\begin{itemize}
\item[a)] $\calE=\calO_C\oplus \calO_C(d)$ with 
$d\in \Pic^0(C)$ non-torsion. Let $C'$ be the section
with $C'\sim C_0-\pi^*d$. Then $K+C_0+C'=0$ and 
$D=d_0 C_0+d'C'$.
\item[b)] $\calE$ is an indecomposable extension of
$\calO_C$ by $\calO_C$, $K+2C_0=0$ and $D=d_0 C_0$.
\end{itemize}
\end{itemize}
\end{thm}

\begin{proof} (of Theorem~\ref{sf}) 
We contract all $(-1)$-curves on which $D$ is numerically 
trivial: we have a birational contraction $f\colon X\to Y$ 
such that $D=f^*(D_Y)$ and $A=K-f^*(K_Y)$ is effective, 
exceptional on $Y$. In particular,
$$
\kappa(X,K+tD)=\kappa(Y,K_Y+tD_Y) \mbox{ for } 
t\in \R.
$$
By construction, $D_Y$ is positive on every 
$K_Y$-negative extremal ray of $Y$. Note that $Y$ is 
not a Del Pezzo surface: otherwise $D_Y$ is semi-ample, 
hence good, by the Base Point Free Theorem. 
Therefore $-K_Y\cdot R\le 1$ for every $K_Y$-negative 
extremal ray $R$ of $Y$. Moreover, $D_Y\cdot R \ge 1$
since $D_Y$ is Cartier. Therefore 
$K_Y+tD_Y$ is nef for $t\ge 2$. 
In particular,
$$
(K_Y+tD_Y)^2=K_Y^2+2(K_Y\cdot D_Y) t+(D_Y^2) t^2 \ge 
0 \mbox{ for } t\ge 2.
$$ 
Therefore either $(K_Y+tD_Y)^2>0$ for $t>2$ (case (1)), 
or $K_Y^2=K_Y\cdot D_Y=D_Y^2=0$. Assume the latter
holds.
By Theorem~\ref{srk}.(4), $D_Y$ is algebraically equivalent
to an effective divisor $D'$. The pairs $(Y,D')$ are 
classified by Theorem~\ref{sk}. Exactly one of the 
following holds:
\begin{itemize}
\item[(i)] $Y$ is a rational surface and there exists
$m\in \N$ such that $K_Y+\frac{1}{m}D_Y\equiv 0$.
\item[(ii)] $Y={\mathbb P}_C(\calE)$, where $C$ is an 
elliptic curve and $\deg(\calE)=0$, and $K_Y+tD_Y\equiv 
0$ for some $0< t\le 2$.
\end{itemize}
\end{proof}


\end{document}